\documentclass[a4paper,12pt]{article}

\usepackage{amsmath,amsfonts,amssymb,amsthm,mathtools,mathrsfs}
\usepackage[cp1251]{inputenc}
\usepackage[T2A]{fontenc}
\usepackage[english]{babel}
\usepackage{hyperref}
\usepackage{dsfont}
\usepackage[all,cmtip]{xy}

\newtheorem{theorem}{Theorem}[section]
\newtheorem{lemma}{Lemma}[section]
\newtheorem{corollary}{Corollary}[section]
\newtheorem{proposition}{Proposition}[section]

\theoremstyle{definition}
\newtheorem{definition}{Definition}[section]
\newtheorem{example}{Example}[section]
\newtheorem{remark}{Remark}[section]

\numberwithin{equation}{section}

\title{Reshetnyak-class mappings and composition operators}
\author{Stepan V. Pavlov\footnote{Department of Mechanics and Mathematics, Novosibirsk State University, 1 Pirogov St, Novosibirsk, 630090, Russian Federation. E-mail: s.pavlov4254@gmail.com},\quad Sergey K. Vodopyanov\footnote{Sobolev Institute of Mathematics,
Siberian Branch of the Russian Academy
of Sciences, 4 Akademika Koptyuga Pr, Novosibirsk, 630090, Russian Federation. E-mail: vodopis@math.nsc.ru}}
\date{}
\begin{document}
\maketitle 
\begin{abstract}
For the Reshetnyak-class homeomorphisms
$\varphi:\Omega\to Y$,
where~$\Omega$
is a~domain in some Carnot group
and~$Y$ is a~metric space,
we obtain an~equivalent description as the mappings
which induce the bounded composition operator
$$
\varphi^*:{\rm Lip}(Y)\to L_q^1(\Omega),
$$
where
$1\leq q\leq \infty$,
as
$\varphi^*u=u\circ\varphi$
for
$u\in{\rm Lip}(Y)$.
We demonstrate the utility of our
approach
by characterizing
the homeomorphisms
$\varphi:\Omega\to\Omega'$
of domains in some Carnot group~$\mathbb G$
which induce the bounded composition operator
$$
\varphi^*: L^1_p(\Omega')\cap {\rm Lip}_{{\rm loc}}(\Omega')\to  L^1_q
(\Omega),\quad 1\leq q \leq p\leq \infty,
$$
of homogeneous Sobolev spaces.
The new proof is much
shorter than the
one already available,
requires a~minimum of tools,
and enables us to obtain
new properties of the homeomorphisms in question.
\end{abstract}

\tableofcontents

\section{Introduction}

In \cite{VodUkhlov1998, VodEvs2022},
a~description is obtained for
the homeomorphisms
$\varphi:\Omega\to\Omega'$
of domains in an~arbitrary Carnot group\footnote{
For the definition of a~Carnot group,
see Section~\ref{prelim-carnot}.
The Euclidean space
$\mathbb{R}^n$
is a~particular case of Carnot groups.}
$\mathbb G$
inducing the bounded composition operator
\begin{equation}\label{primres}
\varphi^*: L^1_p(\Omega')\cap {\rm Lip}_{{\rm loc}}(\Omega')\to  L^1_q
(\Omega),\quad 1\leq q \leq p<\infty,
\end{equation}
as the change-of-variable rule:
$\varphi^*u=u\circ\varphi$
for
$u\in  L^1_p(\Omega')\cap {\rm Lip}_{{\rm loc}}(\Omega')$.

\begin{theorem}\label{descrCarnot0}
The operator~\eqref{primres}
is bounded
if and only if
\begin{itemize}
\item[\rm(1)]
$\varphi\in W^1_{q,{\rm loc}}(\Omega;\mathbb{G});$
\item[\rm(2)]
$D_h\varphi=0$
a.e.\ on
$\{x\in\Omega\mid \det\widehat{D}\varphi(x)=0\}$;
\item[\rm(3)]
$K_{p}(\cdot,\varphi)\in L_{\sigma}(\Omega)$,
where
$\frac{1}{\sigma}=\frac{1}{q}-\frac{1}{p}$,
as well as
$\sigma=\infty$
for
$q=p$,
while the outer distortion function
$K_{p}(\cdot,\varphi)$
for
$p\in[1;\infty)$
is defined as
$$
\Omega\ni x \mapsto K_{p}(x,\varphi) =
\begin{cases}
\dfrac{|D_h \varphi(x)|}
{|\det \widehat{D}\varphi(x)|^{\frac{1}{p}}}
& \textit{if } \det \widehat{D} \varphi(x) \neq 0, \\
0
& \textit{otherwise}.
\end{cases}
$$
\end{itemize}
Furthermore,
$$
\alpha_{q,p}\|K_p(\cdot,\varphi)\mid L_\sigma(\Omega)\|\leq \|\varphi^*\|\leq \|K_p(\cdot,\varphi)\mid L_\sigma(\Omega)\|,
$$
where
$\alpha_{q,p}\in (0;1)$
is a~constant.
\end{theorem}

In Section~\ref{p<oo}
(Theorem~\ref{descrCarnot})
we present a~new proof that
conditions (1)--(3) are necessary
and refine the estimate for the composition operator:
instead of the two-sided inequality
we establish that
$$
\|\varphi^*\|=\|K_p(\cdot,\varphi)\mid L_\sigma(\Omega)\|.
$$
In the case of the Euclidean space
$\mathbb{G}=\mathbb R^n$
the equality between the norms of the composition operator
and the distortion function
is obtained in~\cite{Vod2022}.

Note that
Theorem \ref{descrCarnot0} leaves out the case
$1\leq q< p=\infty$.
It is discussed in the first half of this article;
moreover,
instead of the composition operator
$$
\varphi^*: L^1_\infty(\Omega')\cap{\rm Lip}_{{\rm loc}}(\Omega')\to  L^1_q(\Omega), \quad 1\leq q <\infty,
$$
we study more general composition operators\footnote{
In the case
$\Omega\subset \mathbb R^n$
the problem of describing the homeomorphisms
$\varphi:\Omega \to Y$
for which the composition operator~\eqref{prim} is bounded
is stated in~\cite[Remark~6]{Vodop2012}.
}
of the form
\begin{equation}\label{prim}
 \varphi^*: {\rm Lip}(Y)\to L_q^1(\Omega),\quad 1\leq q<\infty,
\end{equation}
where~$\Omega$
is a~domain in some Carnot group~$\mathbb{G}$
and
$(Y,d)$
is a~metric space,
while
$\varphi:\Omega \to Y$
is a~homeomorphism.

Our main result is Theorem~\ref{descr},
which links 
the operators~\eqref{prim}
and Reshetnyak-class mappings,
see Definition~\ref{R1p}.
Its proof rests on the new Lemma~\ref{mainIneq},
stated and proved by S.~V.~Pavlov.

An~important role in these questions
is played by the construction~\cite{VodUkhlov1998}
of the quasi-additive set function associated to the composition operator.

\section{Preliminaries}\label{prelim}
\subsection{Carnot groups}\label{prelim-carnot}
Recall that
a~\emph{stratified graded nilpotent group},
or a~\emph{Carnot group},
see \cite[Ch.~1]{FS} for instance,
is a~simply-connected connected Lie group~$\mathbb{G}$
whose Lie algebra~$\mathfrak{g}$
of left-invariant vector fields
decomposes as a~direct sum
$\mathfrak{g}
= \mathfrak{g}_1 \oplus \mathfrak{g}_2 \oplus \dots \oplus \mathfrak{g}_m$
of subspaces
$\mathfrak{g}_i$
satisfying
$[\mathfrak{g}_1, \mathfrak{g}_i] = \mathfrak{g}_{i+1}$
for
$i=1,\dots,m-1$
and
$[\mathfrak{g}_1, \mathfrak{g}_m] = \{ 0 \}$.

Fix some inner product in~$\mathfrak{g}$.
The subspace
$\mathfrak{g}_1\subset \mathfrak{g}$
is called the \emph{horizontal space} of~$\mathfrak{g}$,
and its elements are
\emph{horizontal vector fields}.
Put
$N = \dim \mathfrak{g}$
and
$n_i = \dim \mathfrak{g}_i$
for
$i = 1, \ldots, m$.
For convenience,
also put
$n = n_1$.
Fix an~orthonormal basis
$X_{i1}, \ldots, X_{in_i}$
of
$\mathfrak{g}_i$.
Since the exponential mapping
$
g = \exp \Big(
\sum\limits_{i=1}^m \sum\limits_{j=1}^{n_i} x_{ij} X_{ij}
\Big)(e),
$
where~$e$
is the neutral element of~$\mathbb{G}$,
is a~global diffeomorphism of~$\mathfrak{g}$
onto~$\mathbb{G}$
\cite[Prop.~1.2]{FS},
we can identify the point
$g \in \mathbb{G}$
with the point
$x = (x_{ij}) \in \mathbb{R}^N$.
Then
$e = 0$
and
$x^{-1} = -x$.
The dilations
$\delta_\lambda$
defined as
$\delta_\lambda (x_{ij}) = (\lambda^i x_{ij})$
are group automorphisms for all
$\lambda > 0$.

A~curve
$\gamma:[a;b]\to\mathbb{G}$
which is absolutely continuous in the Riemannian metric
is called {\it horizontal}
whenever
$\dot{\gamma} (t)\in \mathfrak{g}_1(\gamma(t))$
for a.e.~$t$.
Refer as the {\it Carnot--Carath\'eodory distance}
$d_{cc}(x,y)$
between two points
$x,\, y\in\mathbb{G}$
to the greatest lower bound  of the lengths
$\int_a^b|\dot{\gamma}(t)|\, dt$
of horizontal curves with endpoints~$x$ and~$y$.
According to the Rashevski\u\i{}--Chow theorem,
see \cite[Section~0.4,\,Section~1.1]{Gromov} for instance,
two arbitrary points can be connected by
a~piecewise smooth horizontal curve of finite length.
The Carnot--Carath\'eo\-dory distance
$d_{cc}$
is a~metric on~$\mathbb{G}$
inducing on it
the original topology,
but not equivalent to the Riemannian metric.

The Lebesgue measure
$dx$
on
$\mathbb{R}^N$
is
a~bi-invariant Haar measure on~$\mathbb{G}$,
and
$d(\delta_\lambda x) = \lambda^\nu dx$,
where
$\nu = \sum\limits_{i=1}^m i\,n_i$
is the~\emph{homogeneous dimension} of the group~$\mathbb{G}$.
The measure is normalized so that
its value on the unit ball
$B(0,1)$
equals~1.
Here
$B(x,r)=\{y\in\mathbb{G}\mid d_{cc}(x,y)<r\}$
is a~ball in the Carnot--Carath\'eodory metric.

\begin{example}
The Heisenberg group
$\mathbb{H}^k = (\mathbb{R}^{2k+1}, *)$
with the group operation
$$
(x, y, z) * (x', y', z')
= \big(
x + x',
y + y',
z + z' + \tfrac{x \cdot y' - x' \cdot y}{2}
\big),
\quad
x, x', y, y' \in \mathbb{R}^k,
\:\:
z, z' \in \mathbb{R},
$$
is a~classical example of a~nonabelian Carnot group.
Its Lie algebra
$\mathfrak{h}^k=\mathfrak{h}^k_1\oplus \mathfrak{h}^k_2$
is formed by the vector fields
$$
X_i = \frac{\partial}{\partial x_i}
- \frac{y_i}{2} \frac{\partial}{\partial z},
\quad
Y_i = \frac{\partial}{\partial y_i}
+ \frac{x_i}{2} \frac{\partial}{\partial z},
\quad
i = 1, \ldots, k,
\quad
Z = \frac{\partial}{\partial z}.
$$
Here
$\mathfrak{h}_1^k = \mathrm{span} \{ X_i, Y_i \mid i = 1, \ldots, k \}$
and
$\mathfrak{h}_2^k = \mathrm{span} \{ Z \}$,
while the only nontrivial Lie brackets are
$[X_i, Y_i] = Z$
for
$i = 1, \ldots, k$.
The homogeneous dimension of
$\mathbb{H}^k$
equals
$\nu = 2k + 2$.
\end{example}

\subsection{Sobolev functions and Reshetnyak classes of
  mappings with values in a~metric space}
\label{prelim-sobol}
Consider a~domain
$\Omega \subset \mathbb{G}$,
which just means a~nonempty connected open set in~$\mathbb{G}$.
For
$p \in [1;\infty)$
the space
$L_p(\Omega)$
consists of all measurable functions
$u : \Omega \to \mathbb{R}$
which are integrable to power~$p$.
The norm on
$L_p(\Omega)$
is defined as
$$
\| u \mid L_p(\Omega)\|= \biggl( \int\limits_\Omega
|u(x)|^p \, dx \biggr)^{\frac{1}{p}}.
$$
The space
$L_\infty(\Omega)$
consists of all measurable essentially bounded functions
$u:\Omega\to\mathbb{R}$.
The norm on
$L_\infty(\Omega)$
is defined in the standard fashion as
$$
\| u \mid L_\infty(\Omega)\|=\operatorname*{ess\,sup}\limits_{x\in\Omega}|u(x)|,
$$
where
$\operatorname*{ess\,sup}\limits_{x\in\Omega}|u(x)|$
is the essential supremum of~$u$
on~$\Omega$.
A~function~$u$
belongs to
$L_{p,\mathrm{loc}}(\Omega)$
for
$p\in [1;\infty]$
whenever
$u \in L_p(K)$
for every compact set
$K \subset \Omega$.

Denote a~basis of the horizontal space
$\mathfrak{g}_1$
by
$X_1=X_{11}, \ldots, X_{n}=X_{1n}$.
Denote by
$\Pi_j$
the hyperplanes
$\{x\in \mathbb G\mid x_j = 0\}$
for
$j = 1, \ldots, n$,
where
$x_j=x_{1j}$
is the horizontal coordinate of the point
$x=(x_{ij})$.
The contraction of
$X_j$
with the volume form
determines the measure
$d\mu_j = \imath(X_j) dx$
on
$\Pi_j$.
Each
$x \in \mathbb{G}$
corresponds to some integral line
$\mathbb{R}\ni t\mapsto \exp ( t X_j)(x)$
intersecting
$\Pi_j$
in a~unique point.
Therefore,
we have a~projection
${\rm Pr}_j:\mathbb{G}\to \Pi_j$.
For $x\in\Pi_j$,
denote by
$$
A_x=\{t\in \mathbb{R}\mid \exp (tX_j)(x)\in A\}
$$
the~$x$-section
of
$A\subset\mathbb{G}$.

Given a~metric space
$(Y,d)$
and a~domain
$\Omega\subset\mathbb{G}$,
a~mapping
$\varphi:\Omega\to Y$
is called {\it measurable}
whenever the preimage
$\varphi^{-1}(M)$
of every Borel set
$M\subset Y$
is Lebesgue measurable.

\begin{definition}\label{ACL-def}
Say that
a~measurable mapping
$\varphi:\Omega\to Y$
is of class
${\rm ACL}(\Omega;Y)$
and call it {\it absolutely continuous on almost all lines}
whenever for each
$j=1,\ldots,n$
we can redefine~$\varphi$
on some measure zero set so that
for
$\mu_j$-a.e.
$x\in {\rm Pr}_j(\Omega)$
the mapping
$$
\Omega_x\ni t\mapsto \varphi(\exp(tX_j)(x))\in Y
$$
is absolutely continuous on every segment\footnote{
Observe that
$d_{cc}(\exp(tX_j)(x),\exp(sX_j)(x))=|t-s|$.}
$[a;b]\subset \Omega_x$.
Put
$\operatorname{ACL}(\Omega)=\operatorname{ACL}(\Omega;\mathbb{R})$.
\end{definition}

It is known\footnote{
See \cite[Theorem 2.7.6]{BBI} for instance,
where this assertion is proved for Lipschitz curves.
The case of absolutely continuous curves and
${\rm ACL}$-mappings
requires only minor modifications.}
that every mapping
$\varphi\in{\rm ACL}(\Omega;Y)$
a.e.\ has {\it metric partial derivatives}
$$
{\rm m} X_j\varphi(x)=\lim\limits_{h\to 0}\frac{d(\varphi(\exp(h X_j)(x)),\varphi(x))}{|h|}=\lim\limits_{h\to +0}\frac{\ell_{\varphi,j}(x,h)}{h},
$$
where
$\ell_{\varphi,j}(x,h)$
is the length of the curve
$t\mapsto \varphi(\exp(tX_j)(x))$
on the interval
$[0;h]$.

For
$p \in [1;\infty]$
the space
$L^{1}_{p}(\Omega)$
consists of all functions
$u \in L_{1,\mathrm{loc}}(\Omega) \cap {\rm ACL}(\Omega)$
whose classical derivatives\footnote{
More exactly,
the derivatives of a~representative for~$u$
which is absolutely continuous on a.e.\ integral line of the fields
$X_1,\ldots, X_n$.
The classical derivatives of this representative
in the directions of the fields
$X_j$
exist a.e.}
$X_j u$
for
$j = 1, \ldots, n$
belong to
$L_p(\Omega)$.
The seminorm of
$u\in L^{1}_{p}(\Omega)$
equals
$$
\| u\mid {L^{1}_{p}(\Omega)}\|
= \big\| \, | \nabla\!_\mathit{h} u | \mid L_p(\Omega)\big\|,
$$
where
$\nabla\!_\mathit{h} u = (X_1 u, \ldots, X_n u)=\sum\limits_{j=1}^n (X_ju) X_j$
is the \emph{horizontal gradient} of~$u$.
Henceforth,
instead of
$\big\| \, | \nabla\!_\mathit{h} u | \mid L_p(\Omega)\big\|$
we write
$\| \nabla\!_\mathit{h}u\mid {L_{p}(\Omega)}\|$.

An~equivalent definition of the space
$L^1_p(\Omega)$
is based on the concept of generalized derivative in the sense of Sobolev:
a~locally summable function
$u_j:\Omega\to\mathbb{R}$
is called the {\it generalized derivative of some function
$u\in L_{1,\mathrm{loc}}(\Omega)$
along of the vector field}
$X_j$,
for
$j = 1,\ldots, n$,
whenever
$$
\int\limits_\Omega u_j(x) v(x)\, dx = - \int\limits_\Omega u(x)X_jv(x)\,dx
$$
for every test function
$v\in C_0^\infty(\Omega)$.
A~locally summable function
$u:\Omega\to\mathbb{R}$
belongs to
$L^1_p(\Omega)$
if and only if
it has the generalized derivatives
$u_j\in L_p(\Omega)$
for
$j=1,\ldots,n$.
Moreover,
$u_j=X_ju$
a.e.,
where
$X_ju$
are the classical derivatives of the function
$u\in \operatorname{ACL}(\Omega)$.

The Sobolev space
$W^1_p(\Omega)$
consists of all functions
$u \in L_p(\Omega) \cap L^1_p(\Omega)$
and is equipped with the norm
$$
\| u \mid W^1_p(\Omega)\|
= \| u \mid L_p(\Omega)\| + \| u \mid L^1_p(\Omega)\|.
$$

Consider a~metric space
$(Y,d)$
and some domain
$\Omega\subset\mathbb{G}$.
For
$1\leq p\leq \infty$
the class
$L_{p}(\Omega;Y)$
consists of all measurable mappings\footnote{
More exactly,
of the classes of such mappings identified in the case of coincidence a.e.}
$\varphi:\Omega\to Y$
satisfying
$$
\| d(\varphi(\cdot),z)\mid L_p(\Omega)\|<\infty
$$
at every point
$z\in Y$.
For
$|\Omega|=\infty$
and
$1\leq p<\infty$
the definition of the class
$L_p(\Omega;Y)$
may depend on the point~$z$,
and so it is meaningless
unless a~``central'' point is chosen in~$Y$
like the neutral element
$e=0$
of a~Carnot group.
The class
$L_{p,{\rm loc}}(\Omega;Y)$
consists of all measurable mappings
$\varphi:\Omega\to Y$
with
$\varphi\in L_p(U;Y)$
for all compact subdomains
$U\Subset\Omega$.

The space
${\rm Lip}(Y)$
consists of all Lipschitz functions
$u:Y\to\mathbb{R}$
and is equipped with the seminorm
$$
{\rm Lip}(u)=\sup\limits_{y_1\neq y_2} \frac{|u(y_1)-u(y_2)|}{d(y_1,y_2)}.
$$

\begin{definition}[{\cite{Resh97},\,\cite{Vodopyanov99}}]\label{R1p}
Say that
a~measurable mapping
$\varphi:\Omega\to Y$
is of {\it Reshetnyak class}
$L^1_p(\Omega;Y)$,
where
$1\leq p\leq \infty$,
whenever
\begin{itemize}
\item[\rm(1)]
$u\circ\varphi\in L^1_p(\Omega)$
for all
$u\in {\rm Lip}(Y)$
and
\item[\rm(2)]
there exists a~function
$g\in L_p(\Omega)$
such that
\begin{equation}\label{UpGrad}
|\nabla_h(u\circ \varphi)|\leq {\rm Lip}(u)g
\end{equation}
a.e.\ in~$\Omega$
for every function
$u\in {\rm Lip}(Y)$.
\end{itemize}
\end{definition}

The class
$W^1_p(\Omega;Y)$
consists of mappings lying in
$L_p(\Omega;Y)\cap L^1_p(\Omega;Y)$.
A~mapping~$\varphi$
lies in
$W^1_{p,{\rm loc}}(\Omega;Y)$
whenever
$\varphi\in W^1_p(U;Y)$
for every compact subdomain
$U\Subset \Omega$.

\begin{remark}
The class
$W^1_p(\Omega;Y)$
was considered  for the first time in~\cite{Resh97}
in the case that~$\Omega$
is a~domain in the Euclidean space~$\mathbb{R}^n$.
\end{remark}

\begin{remark}
If~$Y$
is the Euclidean space
$\mathbb{R}^k$
then the membership of a mapping
$\varphi$
in the Reshetnyak class
$L^1_p(\Omega;\mathbb{R}^k)$
or
$W^1_p(\Omega;\mathbb{R}^k)$
is equivalent to the property that
all coordinate functions
$\varphi_1,\ldots,\varphi_k$
of~$\varphi$
lie in
$L^1_p(\Omega)$
or
$W^1_p(\Omega)$.
\end{remark}

\begin{remark}
In the case of a~separable space~$Y$
we can show that
the class
$L^1_p(\Omega;Y)$
does not enlarge
when we impose conditions~(1) and~(2),
instead of all Lipschitz functions,
only on the~$1$-Lipschitz distance functions
$u=u_{z_i}=d(\cdot,z_i)$,
where
$z_i$
runs over a~countable dense subset of~$Y$.
The stock of functions~$g$
satisfying~\eqref{UpGrad}
also remains the same.
The arguments exactly repeat
\cite[Prop.~4.1]{Vodop2000} and \cite[Corollary 1]{Resh97},
where this property is established for
$W^1_p(\Omega;Y)$.
\end{remark}

Condition~(1) clearly implies that
$\varphi\in L_{1,{\rm loc}}(\Omega;Y)$
because the functions
$d(\varphi(\cdot),z)$,
where
$z\in Y$,
are locally summable.
It is known \cite{KVP}, \cite[6.10, Theorem 17]{KantAkil} that
among the functions~$g$
satisfying~\eqref{UpGrad}
there is a~minimal one,
called the {\it upper gradient} of the mapping~$\varphi$
and denoted by
$|\nabla_0\varphi|$.
The minimality of
$|\nabla_0\varphi|$
means that
if
$g\in L_p(\Omega)$
satisfies~\eqref{UpGrad}
then
$|\nabla_0\varphi|\leq g$
a.e.

Consider two Carnot groups~$\mathbb{G}$ and~$\widetilde{\mathbb{G}}$,
as well as a~domain
$\Omega \subset \mathbb{G}$.
If
$\varphi \in {\rm ACL}(\Omega; \widetilde{\mathbb{G}})$
then
$X_j \varphi(x) \in \widetilde{\mathfrak{g}}_1(\varphi(x))$
for a.e.\
$x \in \Omega$
\cite[Prop.~4.1]{Pansu}.
The matrix
$D_\mathit{h} \varphi(x) = (X_i \varphi_j)$,
where
$i = 1, \ldots, n$
and
$j = 1, \ldots, \widetilde{n}$,
determines the linear operator
$D_\mathit{h} \varphi(x) : \mathfrak{g}_1 \to \widetilde{\mathfrak{g}}_1$
called the \emph{horizontal differential} of~$\varphi$.
It is known \cite[Theorem~1.2]{Vodop2000} that
for a.e.\
$x \in \Omega$
the linear operator
$D_\mathit{h}\varphi(x)$
is defined and extends
to a~homomorphism
$\widehat{D} \varphi (x): \mathfrak{g} \to \widetilde{\mathfrak{g}}$
of Lie algebras,
which we can also regard as the linear operator
$\widehat{D} \varphi (x): T_x \mathbb{G} \to T_{\varphi(x)} \widetilde{\mathbb{G}}$.
The norms of these operators satisfy
$$
\big| D_\mathit{h} \varphi(x) \big|
\le \big| \widehat{D} \varphi(x) \big|
\le C \big| D_\mathit{h} \varphi(x) \big|,
$$
where~$C$
depends only on the structure of the groups.
Here the norm of the operator
$\widehat{D} \varphi(x)$
is defined as
$$
\sup\big\{\widetilde{d}_{cc}\big(\widehat{D}\varphi(x)\langle X\rangle\big)\mid X\in\mathfrak{g},\, d_{cc}(X)\leq 1\big\},
$$
where we put
$d_{cc}(X)=d_{cc}(\exp (X),0)$
and
$\widetilde{d}_{cc}(\widetilde{X})=\widetilde{d}_{cc}(\widetilde{\exp}(\widetilde{X}),0)$
for
$X\in\mathfrak{g}$
and
$\widetilde X\in\widetilde{\mathfrak{g}}$
for brevity.
Corresponding to the homomorphism
$\widehat{D} \varphi(x)$,
there is also the group homomorphism
$D_\mathcal{P} \varphi (x)
= \widetilde{\exp} \circ \widehat{D} \varphi(x) \circ \exp^{-1}$
known as the \emph{Pansu differential},
which is the approximative differential of~$\varphi$
with respect to the group structure \cite{Vodop2000}.

Equivalent descriptions for the mappings of Reshetnyak class
$W^1_p(\Omega;Y)$
are established in \cite{Vodopyanov99} and \cite{Vodop2000}.
The necessary arguments carry over to the case of class
$L^1_p(\Omega;Y)$
without any change.

\begin{proposition}[{\cite[Section~5]{Vodopyanov99},\,\cite[Prop.~4.2]{Vodop2000}}]\label{R1p-ACL}
Given a~complete separable metric space
$(Y,d)$,
a~domain~$\Omega$
in some Carnot group
$\mathbb{G}$,
and
$1\leq p\leq\infty$,
the following conditions on a~mapping
$\varphi:\Omega\to Y$
are equivalent:
\begin{itemize}
\item[\rm(1)]
$\varphi\in L^1_p(\Omega;Y)$;
\item[\rm(2)]
$\varphi\in L_{1,{\rm loc}}(\Omega;Y)\cap {\rm ACL}(\Omega;Y)$
and
${\rm m}X_j\varphi\in L_p(\Omega)$
for
$j=1,\ldots,n$.
\end{itemize}

If
$X_1,\ldots,X_n$
is an~orthonormal basis for the space
$\mathfrak{g}_1$
then
${\rm m} X_j\varphi\leq |\nabla_0 \varphi|\leq
\left(\sum\limits_{i=1}^n({\rm m} X_i\varphi)^2\right)^{1/2}$
a.e.

If in addition
$(Y,d)=(\widetilde{\mathbb{G}},\widetilde{d}_{cc})$
is a~Carnot group endowed with the Carnot--Cara\-th\'eodory metric
then each of conditions (1) and (2)
is equivalent to the following:
\begin{itemize}
\item[\rm(3)]
$\varphi\in L_{1,{\rm loc}}(\Omega;\widetilde{\mathbb{G}})\cap {\rm ACL}(\Omega;\widetilde{\mathbb{G}})$
and
$|D_h\varphi|\in L_p(\Omega)$.
\end{itemize}
Furthermore,
${\rm m}X_j\varphi=|X_j\varphi|$
a.e.\ for
$j=1,\ldots,n$.
\end{proposition}

\subsection{A quasi-additive set function}
\label{prelim-quasiadd}
A~set function
$\Phi:\mathcal{O}(Y)\to [0;+\infty]$
defined on the system
$\mathcal{O}(Y)$
of all open subsets of a~metric space
$(Y,d)$
is called {\it quasi-additive}
whenever:
\begin{itemize}
\item[(1)]
for every point
$y\in Y$
there exists
$\delta>0$
such that
$\Phi(B(y,r))<\infty$
for
$r<\delta$,
where
$B(y,r)=\{z\in Y\mid d(z,y)<r\}$
is an~open ball;
\item[(2)]
every finite disjoint collection
$V_1,\ldots, V_k\in\mathcal{O}(Y)$
with
$V_1\cup\ldots\cup V_k\subset V\in\mathcal{O}(Y)$
satisfies
$$
\sum\limits_{j=1}^k\Phi(V_j)\leq \Phi(V).
$$
\end{itemize}

\begin{proposition}[{\cite[Corollary 5]{VU2003}}]\label{quasiadditprop}
If
$\Phi:\mathcal{O}(\Omega)\to [0;+\infty]$
is a~quasi-additive function,
where~$\Omega$
is a~domain in some Carnot group~$\mathbb{G}$,
then:
\begin{itemize}
\item[\rm(a)]
for a.e.
$x\in\Omega$
the derivative
$$
\Phi'(x)=\lim\limits_{\delta\to 0,B_\delta\ni x}\frac{|\Phi(B_\delta)|}{|B_\delta|}
$$
exists and is finite,
where
$B_\delta$
is an~arbitrary radius~$\delta$ ball in the metric
$d_{cc}$
containing the point~$x$;
\item[\rm(b)]
$\Phi'$
is a~measurable function;
\item[\rm(c)]
for every
$U\in \mathcal{O}(\Omega)$
we have
$$
\int\limits_U \Phi'(x)\, dx\leq \Phi(U).
$$
\end{itemize}
\end{proposition}

\begin{example}
Consider a~homeomorphism
$\varphi:\Omega\to\Omega'$
between two domains
$\Omega$,
$\Omega'\subset \mathbb{G}$.
The function of open sets
$$
\mathcal{O}(\Omega)\ni U\mapsto |\varphi(U)|
$$
is (quasi-)additive.
By Proposition~\ref{quasiadditprop},
for a.e.
$x\in\Omega$
the limit
\begin{equation}\label{Jacob}
J(x,\varphi)=\lim\limits_{r\to 0}\frac{|\varphi(B(x,r))|}{|B(x,r)|}
\end{equation}
exists,
and it is called the {\it spatial derivative} of $\varphi$.
If $\varphi$
is of class
$W^1_{1,{\rm loc}}(\Omega;\mathbb{G})$
then its spatial derivative coincides with
the determinant of its Pansu differential:
$J(\cdot,\varphi)=|\det \widehat{D}\varphi|$
a.e.,
see \cite[Lemma 2.1]{VodEvs2022} for instance.
\end{example}

\section{An auxiliary property}
\label{techn}

The proof of the following available property given below
is a~straightforward generalization of the arguments for
$\mathbb{G}=\mathbb{R}^n$
in \cite[Section~4.2, Theorem 4]{EG}.
Given a~function
$u:\Omega\to\mathbb{R}$,
denote by
$\chi_+$
and
$\chi_-$
the characteristic functions of the sets
$\{x\mid u(x)>0\}$
and
$\{x\mid u(x)<0\}$
respectively.

\begin{proposition}\label{chain}
If
$u\in L^1_q(\Omega)$,
where~$\Omega$
is a~domain in some Carnot group~$\mathbb{G}$
and
$1\leq q<\infty$,
while the derivative of some function
$F\in C^1(\mathbb{R})$
is bounded
then the following claims hold.
\begin{itemize}
\item[\rm(1)]
$F\circ u\in L^1_q(\Omega)$;
furthermore,
$\nabla_h(F\circ u)(x)=F'(u(x))\nabla_h u(x)$
for a.e.
$x\in\Omega$.
\item[\rm(2)]
The functions
$u^+=\max\{u,0\}$,
$u^-=-\min\{u,0\}$,
and
$|u|$
lie in
$L^1_q(\Omega)$;
moreover,
$$
\nabla_h u^+=\chi_+ \nabla_h u,\quad \nabla_h u^-
=-\chi_- \nabla_h u,\quad \nabla_h |u|={\rm sgn}\, u\cdot \nabla_h u
$$
a.e.\
This implies that
$\nabla_h u=0$
a.e.\ on 
$\{x\mid u(x)=0\}$,
and that
$|\nabla_h u|=|\nabla_h |u||$
a.e.
\item[\rm(3)]
The cutoff
$$
u_M(x)={\rm cut}_M u(x)=\begin{cases}
{\rm sgn}\, u(x)\cdot M
& \text{if\/ } |u(x)|\geq M,\\
u(x)
& \text{if\/ } |u(x)|<M,
\end{cases}
$$
of the function~$u$
belongs to
$L^1_q(\Omega)$
and
$\nabla_h u_M\to \nabla_h u$
in
$L_q(\Omega)$
as
$M\to+\infty$.
\end{itemize}
\end{proposition}

\begin{proof}
(1)
There exists a~sequence
$\{u_l\in C^\infty(\Omega)\cap L^1_q(\Omega)\}$
converging to~$u$
a.e.\ in
$L_{1,{\rm loc}}(\Omega)$
and in
$L^1_{q}(\Omega)$,
see \cite{Garof} or \cite[Lemma 1]{IsangVod2010} for instance.
For
$\psi\in C_0^\infty(\Omega)$,
since~$F'$
is continuous and bounded,
we have
\begin{multline*}
\int\limits_{\Omega}(F\circ u)(x)X_j\psi(x)\, dx=\lim\limits_{l\to\infty}\int\limits_{\Omega}(F\circ u_l)(x)X_j\psi(x)\, dx\\
=-\lim\limits_{l\to\infty} \int\limits_\Omega F'(u_l(x))(X_ju_l(x))\psi(x)\, dx=-\int\limits_\Omega F'(u(x))(X_ju(x))\psi(x)\, dx.
\end{multline*}
The latter means that
$(F'\circ u)\cdot X_ju\in L_q(\Omega)$
is the generalized derivative of the function
$F\circ u$
along the field~$X_j$.

(2)
Consider the functions
$F^+_\varepsilon(t)=\big((t^2+\varepsilon^2)^{1/2}-\varepsilon\big)\chi_{[0;\infty)}(t)$
for
$\varepsilon>0$.
Since
$F^+_\varepsilon\in C^1(\mathbb{R})$
and
$|(F^+_\varepsilon)'|\leq 1$,
by claim~(1) for
$\psi\in C_0^\infty(\Omega)$
we have
\begin{multline*}
\int\limits_\Omega u^+(x)X_j\psi(x)\, dx=\lim\limits_{\varepsilon\to 0}\int\limits_\Omega F^+_\varepsilon(u(x))X_j\psi(x)\, dx\\
=-\lim\limits_{\varepsilon\to 0}\int\limits_\Omega (F^+_\varepsilon)'(u(x))(X_ju(x))\psi(x)\, dx=-\int\limits_\Omega\chi_+(x)(X_ju(x))\psi(x)\, dx.
\end{multline*}
Hence,
$\nabla_h u^+=\chi_+ \nabla_h u$.
Similarly we can establish that
$\nabla_h u^{-}=-\chi_{-} \nabla_h u$.
The formula for
$\nabla_h|u|$
follows from the equality
$|u|=u^+ +u^-$.
The property that
$\nabla_h u=0$
a.e.\ on
$u^{-1}(0)=\{x\in \Omega\mid u(x)=0\}$
follows from the equalities
$$
\nabla_h u=\nabla_h(u^+- u^-)=\chi_+ \nabla_h u+ \chi_- \nabla_h u=\chi_{\Omega\setminus u^{-1}(0)}\nabla_h u.
$$

(3)
Since
$u_M=((u-M)^- +2M)^+ - M$,
it follows that
$\nabla_h u_M=\chi_{u^{-1}((-M;M))}\nabla_h u$.
Since the sets
$u^{-1}((-M;M))$
monotonely increase,
$\Omega=\bigcup\limits_{M>0}u^{-1}((-M;M))$,
and
$\nabla_h u\in L_q(\Omega)$,
we see that
the integrals
$$
\int\limits_{\Omega}|\nabla_h u_M(x)-\nabla_h u(x)|^q\, dx=\int\limits_{\Omega\setminus u^{-1}((-M;M))}|\nabla_h u(x)|^q\,dx
$$
vanish as
$M\to +\infty$.
\end{proof}

\section{A functional description of Reshetnyak-class homeomorphisms}

As above,
consider a~domain~$\Omega$
in some Carnot group~$\mathbb{G}$
and a~metric space
$(Y,d)$.
If a~mapping
$\varphi:\Omega\to Y$
belongs to
$L^1_q(\Omega;Y)$,
where
$1\leq q\leq  \infty$,
then it induces the bounded composition operator
\begin{equation}\label{CompOp}
\varphi^*:{\rm Lip}(Y)\to L^1_q(\Omega),\quad \varphi^*u=u\circ\varphi,
\end{equation}
where
$L^1_q(\Omega)$
is equipped with the seminorm
$f\mapsto \|\nabla_h f\mid L_q(\Omega)\|$,
while
${\rm Lip}(Y)$
stands for the seminorm
$u\mapsto {\rm Lip}(u)$.
Indeed,
since
$|\nabla_h(u\circ\varphi)|\leq {\rm Lip}(u)|\nabla_0\varphi|\in L_q(\Omega)$
a.e.,
it follows that
$$
\|u\circ\varphi\mid L^1_q(\Omega)\|\leq \|\, |\nabla_0\varphi|\mid L_q(\Omega)\|\cdot {\rm Lip}(u).
$$
From this we deduce an~estimate for the norm of operator \eqref{CompOp}:
$$
\|\varphi^*\|\leq \|\, |\nabla_0\varphi|\mid L_q(\Omega)\|.
$$
Assume the converse:
a~certain mapping
$\varphi:\Omega\to Y$
induces the bounded composition operator~\eqref{CompOp}
with
$1\leq q<\infty$.
Is it true that~$\varphi$
is a~mapping of class
$L^1_q(\Omega;Y)$?
We obtain the positive answer to this question
in the case that~$\varphi$
is a~homeomorphism.

Given an~open set
$V\subset Y$,
put\footnote{
Here
$\operatorname{spt} u=\overline{\{y\in Y\mid u(y)\neq 0\}}$
stands for the support of the function
$u:Y\to\mathbb{R}$.}
\begin{multline}\label{Phi-defi}
\Phi(V)=\sup\Big\{\int_\Omega |\nabla_h(u\circ \varphi)(x)|^q\, dx\,  \big\vert\, u\in {\rm Lip}(Y),\\ {\rm Lip}(u)\leq 1,\, {\rm dist}(\operatorname{spt} u,Y\setminus V)>0\Big\}.
\end{multline}
In Lemma~\ref{quasiaddit} we verify that
the set function~$\Phi$
is quasi-additive.

\begin{remark}
We put
${\rm dist}(E,\varnothing)=+\infty$.
Therefore,
for
$V=Y$
the supremum in the definition of~$\Phi$
is taken over all~$1$-Lipschitz functions~$u$.
Hence,
$\Phi(Y)=\|\varphi^*\|^q<\infty$.
We emphasize once more that~$\Phi$
is monotone with respect to inclusion and
$\Phi(\varnothing)=0$.
\end{remark}

\begin{theorem}\label{descr}
Given a~domain~$\Omega$
in some Carnot group~$\mathbb{G}$
and a~metric space
$(Y,d)$,
if a~homeomorphism
$\varphi:\Omega \to Y$
with
$\varphi(\Omega)=Y$
induces the bounded composition operator
$$
\varphi^*:{\rm Lip}(Y)\to L^1_q(\Omega),\quad \varphi^*u=u\circ\varphi,\quad 1\leq q<\infty
$$
then
$\varphi\in L^1_q(\Omega;Y)$.
Moreover,
$$
\int\limits_U |\nabla_0\varphi|(x)^q\, dx=\Phi(\varphi(U))
$$
for every open set
$U\subset\Omega$.
In~particular\footnote{
Here
$(\Phi\circ\varphi)'$
is the derivative of the quasi-additive function
$\Phi\circ\varphi:U\mapsto \Phi(\varphi(U))$,
see Lemma~\ref{quasiaddit} and Proposition~\ref{quasiadditprop}.},
$
|\nabla_0\varphi|^q=(\Phi\circ\varphi)^\prime
$
a.e.\
and
$$
\|\, |\nabla_0\varphi|\mid L_q(\Omega)\|=\|\varphi^*\|.
$$
\end{theorem}

The new Lemma~\ref{mainIneq} proved below contains a~key argument.
By Definition \eqref{Phi-defi} of the function~$\Phi$,
all functions
$u\in{\rm Lip}(Y)$
with
${\rm dist}(\operatorname{spt} u,Y\setminus V)>0$
satisfy\footnote{
The integral in the left-hand side coincides with
$\|\nabla_h (u\circ\varphi)\mid L_q(\Omega)\|^q$.
Indeed,
$\nabla_h(u\circ\varphi)=0$
a.e.\ outside
$\varphi^{-1}(V)$,
so that
$u\circ\varphi=0$
outside
$\varphi^{-1}(V)$,
see assertion~(2) of Proposition~\ref{chain}.}
\begin{equation}\label{Phii}
\int\limits_{\varphi^{-1}(V)}|\nabla_h(u\circ \varphi)(x)|^q\, dx\leq {\rm Lip}(u)^q\Phi(V).
\end{equation}
Let us show how in this estimate
we can remove the condition on the support of~$u$.

\begin{lemma}\label{meas}
Consider a~domain~$\Omega$
in some Carnot group~$\mathbb{G}$,
a~metric space
$(Y,d)$,
a~mapping
$\varphi:\Omega\to Y$,
and an~open subset\/~$V$ of\/~$Y$.

If for every
$u\in {\rm Lip}(Y)$
with
${\rm dist}(\operatorname{spt} u, Y\setminus V)>0$
the composition
$u\circ\varphi$
is measurable
then
$\varphi^{-1}(V)$
is a~measurable set
and
$\varphi:\varphi^{-1}(V)\to V$
is a~measurable mapping.
\end{lemma}

\begin{proof}
Fix an~open subset
$W\subset V$.
The function
$u_W(y)={\rm dist}(y,Y\setminus W)$
is~$1$-Lipschitz.
Define also the~$1$-Lipschitz function
$u^\varepsilon_W(y)=\max\{u_W(y)-\varepsilon,0\}$.
Since
$u^\varepsilon_W(y)>0$
if and only if
${\rm dist}(y,Y\setminus W)>\varepsilon$,
it follows that
$$
{\rm dist}(\operatorname{spt}u^\varepsilon_W,Y\setminus V)\geq {\rm dist}(\operatorname{spt} u^\varepsilon_W,Y\setminus W)\geq \varepsilon>0.
$$
By assumption,
the functions
$u^\varepsilon_W\circ\varphi$
are measurable;
consequently,
so is the set
$\{x\in\Omega\mid (u^\varepsilon_W\circ\varphi)(x)>0\}$.
The measurability of
$\varphi^{-1}(W)$
follows from the equalities
$$
\varphi^{-1}(W)=\{x\in\Omega\mid (u_W\circ\varphi)(x)>0\}=\bigcup\limits_{n=1}^\infty \{x\in\Omega\mid (u^{1/n}_W\circ\varphi)(x)>0\}.
$$
Since the preimage of each open set is measurable,
we infer that
the preimage of every Borel subset of~$V$
is measurable,
which means that
the mapping
$\varphi\vert_{\varphi^{-1}(V)}$
is measurable.
\end{proof}

\begin{lemma}\label{mainIneq}
Consider a~domain
$\Omega$
in some Carnot group~$\mathbb{G}$,
a~metric space
$(Y,d)$,
a~mapping
$\varphi:\Omega\to Y$,
an~open subset
$V\subset Y$,
and a~number
$A\geq 0$.

Take
$1\leq q<\infty$
and suppose that
every~$1$-Lipschitz function
$u:Y\to\mathbb{R}$
with
${\rm dist}(\operatorname{spt} u,Y\setminus V)>0$
satisfies
$u\circ \varphi\in L^1_q(\Omega)$
and
\begin{equation}\label{B<2B}
\int\limits_{\varphi^{-1}(V)}|\nabla_h (u\circ\varphi)(x)|^q\, dx\leq A.
\end{equation}
Then
$\eqref{B<2B}$
for every
$1$-Lipschitz function
$u:Y\to\mathbb{R}$
with
$u\circ\varphi\in L^1_q(\Omega)$.
\end{lemma}

\begin{proof}
{\bf Step 1.}
Consider an~arbitrary~$1$-Lipschitz function
$u:Y\to\mathbb{R}$
with
$u\circ\varphi\in L^1_q(\Omega)$.
Assume firstly that~$u$
is bounded on~$V$.
Put
$M=\sup\limits_{y\in V} |u(y)|<\infty$.
Using {\it symmetrization with respect to a~level curve},
construct a~$1$-Lipschitz function
$\overline{u}$
with the following properties:
$\overline{u}\circ\varphi\in L^1_q(\Omega)$
with
$|\nabla_h (u\circ \varphi)|=|\nabla_h (\overline{u}\circ\varphi)|$
a.e.,
and
$|\overline{u}|\leq \frac{M}{2}$
on~$V$.
To this end,
consider the piecewise linear function
$$
{\rm sym}_{M}(t)=
\begin{cases}
-t-M
& \text{if } t<-\frac{M}{2},\\
t
& \text{if } -\frac{M}{2}\leq t\leq \frac{M}{2},\\
M-t
& \text{if } t>\frac{M}{2}
\end{cases}
$$
of real argument
$t\in [-M;M]$.
Since the symmetrizing function can be expressed as the composition
${\rm sym}_M(t)=\big|M-|t-\frac{M}{2}|\big|-\frac{M}{2}$
of absolute value functions and translations,
where
$|t|\leq M$,
by assertion~(2) of Proposition~\ref{chain}
the function
$\overline{u}={\rm sym}_M\circ u$
satisfies
$|\nabla_h (u\circ \varphi)|=|\nabla_h (\overline{u}\circ\varphi)|$
a.e.\
because
$\overline{u}\circ\varphi={\rm sym}_M\circ (u\circ\varphi)$,
where
$u\circ\varphi\in L^1_q(\Omega)$
by assumption.
The inequality
$|\overline{u}|\leq \frac{M}{2}$
on~$V$
follows from
$|{\rm sym}_M(t)|\leq \frac{M}{2}$
for
$|t|\leq M$.

Repeating this process sufficiently many times,
we obtain a~$1$-Lipschitz function~$\overline{u}$
on~$Y$
with
$\overline{u}\circ\varphi\in L^1_q(\Omega)$
and
$|\nabla_h (u\circ \varphi)|=|\nabla_h (\overline{u}\circ\varphi)|$
a.e.,
while
$|\overline{u}|\leq \delta$
on~$V$,
where~$\delta$
is an~arbitrary prescribed positive number.
Refer to~$\overline{u}$
as a~{\it Lipschitz refinement}
of~$u$.

{\bf Step 2.}
Proceed to justify~
\eqref{B<2B}
for every~$1$-Lipschitz function~$u$
bounded on~$V$
and satisfying
$u\circ\varphi\in L^1_q(\Omega)$.
The family
$$
K_\delta=\{y\in Y\mid {\rm dist}(y,Y\setminus V)>\delta\},\ \delta>0,
$$
of open sets
is monotone,
while its union
$\bigcup\limits_{\delta >0}K_\delta$
equals~$V$.
Therefore,
for every
$\varepsilon>0$
there is a~sufficiently small
$\delta>0$
such that\footnote{
Henceforth we use the measurability of the sets
$\varphi^{-1}(K_\delta)$
obtained in Lemma~\ref{meas}.}
$$
\int\limits_{\varphi^{-1}(K_{2\delta})}|\nabla_h (u\circ\varphi)(x)|^q\, dx
\geq \int\limits_{\varphi^{-1}(V)}|\nabla_h (u\circ\varphi)(x)|^q\, dx-\varepsilon.
$$
As Step~1 shows,
there exists a~$1$-Lipschitz function
$\overline{u}:Y\to\mathbb{R}$
with
$|\overline{u}|\leq \delta$
on~$V$,
while
$\overline{u}\circ\varphi\in L^1_q(\Omega)$
and
$|\nabla_h (u\circ \varphi)|=|\nabla_h (\overline{u}\circ\varphi)|$
a.e.
Let us construct a~cutoff
$\overline{u}_1$
of the refinement~$\overline{u}$.
Put
$\overline{u}_1=\overline{u}$
on
$K_{2\delta}$
and
$\overline{u}_1=0$
on
$Y\setminus K_\delta$.
Since
$|\overline{u}_1|\leq \delta$
on
$K_{2\delta}$
and\footnote{
${\rm dist}(K_{2\delta}, Y\setminus K_\delta)=\infty$
if
$K_\delta=Y$.}
${\rm dist}(K_{2\delta}, Y\setminus K_\delta)\geq \delta$,
it follows that
the function
$\overline{u}_1$
defined outside the annulus
$K_{\delta}\setminus K_{2\delta}$
is~$1$-Lipschitz.
Indeed,
if
$y_1\in K_{2\delta}$
and
$y_2\in Y\setminus K_\delta$
then
$$
|\overline{u}_1(y_1)-\overline{u}_1(y_2)|=|\overline{u}(y_1)|\leq \delta\leq{\rm dist}(K_{2\delta}, Y\setminus K_\delta)\leq d(y_1,y_2).
$$
In the cases that
$y_1$,
$y_2\in K_{2\delta}$
and
$y_1$,
$y_2\in Y\setminus K_\delta$,
a~similar estimate for the increments of
$\overline{u}_1$
is obvious.

By Kirszbraun's theorem\footnote{
See \cite[Section~3.1.1, Theorem 1]{EG} for instance.}
$\overline{u}_1$
extends to a~$1$-Lipschitz function defined on the whole of~$Y$.
Denote this continuation by the same symbol.
By construction the support of
$\overline{u}_1$
is included into
$\overline{K_\delta}$,
and so
${\rm dist}(\operatorname{spt} \overline{u}_1,Y\setminus V)
\geq {\rm dist}(\overline{K_\delta},Y\setminus V)\geq \delta>0$.
The latter property of the support of
$\overline{u}_1$
and the assumption of the lemma show that
the composition
$\overline{u}_1\circ\varphi$
belongs to
$L^1_q(\Omega)$,
and
$\overline{u}_1$
satisfies 
\eqref{B<2B}.
We have
\begin{multline*}
\int\limits_{\varphi^{-1}(V)}|\nabla_h (u\circ\varphi)(x)|^q\, dx-\varepsilon\leq \int\limits_{\varphi^{-1}(K_{2\delta})}|\nabla_h (u\circ\varphi)(x)|^q\, dx\\
=\int\limits_{\varphi^{-1}(K_{2\delta})}|\nabla_h(\overline{u}\circ\varphi)(x)|^q\, dx=\int\limits_{\varphi^{-1}(K_{2\delta})}|\nabla_h (\overline{u}_1\circ\varphi)(x)|^q\, dx\\
\leq \int\limits_{\varphi^{-1}(V)}|\nabla_h (\overline{u}_1\circ\varphi)(x)|^q\, dx\leq A.
\end{multline*}
Here we use the property that
$\nabla_h(\overline{u}_1\circ \varphi)=\nabla_h(\overline{u}\circ \varphi)$
a.e.\ on the set
$$
\{x\in \Omega\mid (\overline{u}_1\circ \varphi)(x)=(\overline{u}\circ \varphi)(x)\}\supset \varphi^{-1}(K_{2\delta}),
$$
see assertion~(2) of Proposition~\ref{chain}.
Since
$\varepsilon>0$
is arbitrary,
we obtain required inequality.

{\bf Step 3.}
Consider an~unbounded
$1$-Lipschitz function
$u:Y\to\mathbb{R}$
with
$u\circ\varphi\in L^1_q(\Omega)$.
Its cutoff
$u_M={\rm cut}_M u={\rm sgn}\, u\cdot \min\{|u|,M\}$
has the following properties:
${\rm Lip}(u_M)\leq 1$
and
$|u_M|\leq M$
on~$Y$,
while
$u_M\circ\varphi=(u\circ\varphi)_M\in L^1_q(\Omega)$.
Step~2 yields
$$
\int\limits_{\varphi^{-1}(V)}|\nabla_h(u_M\circ\varphi)(x)|^q\, dx\leq A,
$$
where the left-hand side tends to
$\int_{\varphi^{-1}(V)}|\nabla_h(u\circ\varphi)(x)|^q\, dx$
as
$M\to+\infty$
by assertion~(3) of Proposition~\ref{chain}
because
$u_M\circ\varphi=(u\circ\varphi)_M$.
\end{proof}

\begin{lemma}\label{quasiaddit}
Given a~domain~$\Omega$
in some Carnot group~$\mathbb{G}$
and a~metric space
$(Y,d)$,
if a~mapping
$\varphi:\Omega\to Y$
induces the bounded composition operator
$$
\varphi^*:{\rm Lip}(Y)\to L^1_q(\Omega),\ \varphi^*u=u\circ\varphi,
$$
where
$1\leq q<\infty$,
then:
\begin{itemize}
\item[\rm(a)]
every open subset
$V\subset Y$
and every
$u\in{\rm Lip}(Y)$
satisfy
$$
\int\limits_{\varphi^{-1}(V)}|\nabla_h (u\circ\varphi)(x)|^q\, dx\leq {\rm Lip}(u)^q \Phi(V),
$$
where the function~$\Phi$
is defined in {\rm\eqref{Phi-defi}};
\item[\rm(b)]
the function~$\Phi$
is quasi-additive.
\end{itemize}
\end{lemma}

\begin{proof}
Claim~(a) follows directly from 
\eqref{Phii},
which is valid for all functions
$u\in {\rm Lip}(Y)$
with
${\rm dist}(\operatorname{spt} u,Y\setminus V)>0$,
and Lemma~\ref{mainIneq} with
$A=\Phi(V)$.

Let us verify that
$\Phi$
is quasi-additive.
Take two disjoint nonempty open subsets
$V_1$
and
$V_2$
of~$Y$.
Consider two
$1$-Lipschitz functions
$u_i$,
for
$i=1,2$,
such that
$$
\int\limits_{\varphi^{-1}(V_i)}|\nabla_h(u_i\circ\varphi)(x)|^q\, dx\geq \Phi(V_i)-\varepsilon,\quad {\rm dist}(\operatorname{spt} u_i,Y\setminus V_i)>0.
$$
For sufficiently large
$M>0$
we have,
see Proposition~\ref{chain},
$$
\int\limits_{\varphi^{-1}(V_i)}|\nabla_h(\tilde{u}_i\circ\varphi)(x)|^q\, dx\geq \int\limits_{\varphi^{-1}(V_i)}|\nabla(u_i\circ\varphi)(x)|^q\, dx-\varepsilon,
$$
where
$\tilde{u}_i={\rm cut}_M u_i$.
The distance
$r={\rm dist}(\operatorname{spt} u_1,\operatorname{spt} u_2)$
is positive.
Consider a~refinement
$\overline{u}_i$
of
$\tilde{u}_i$
such that
$|\overline{u}_i|\leq r/2$
on~$Y$,
while
${\rm Lip}(\overline{u}_i)\leq 1$
and
$|\nabla_h(\overline{u}_i\circ\varphi)|=|\nabla_h(\tilde{u}_i\circ\varphi)|$
a.e.\ on~$\Omega$;
see Step~1 of the proof of Lemma~\ref{mainIneq}.
The function
$u=\overline{u}_1+\overline{u}_2$
coincides with
$\overline{u}_i$
on
$\operatorname{spt} \overline{u}_i\subset\operatorname{spt}\tilde{u}_i=\operatorname{spt}\, u_i$,
for
$i=1,2$,
and vanishes at the remaining points of~$Y$.
The function~$u$
is~$1$-Lipschitz
by the choice of the refinement level:
if
$y_i\in \operatorname{spt} u_i$
for
$i=1,2$
then
$$
|u(y_1)-u(y_2)|=|\overline{u}_1(y_1)-\overline{u}_2(y_2)|\leq \frac{r}{2}+\frac{r}{2}={\rm dist}(\operatorname{spt} u_1,\operatorname{spt} u_2)\leq d(y_1,y_2),
$$
while for other locations of the pairs of points
$y_i$
the same estimate for the increment of~$u$
is obvious.
It is clear that
${\rm dist}(\operatorname{spt} u,Y\setminus(V_1\cup V_2))>0$.
Now the quasi-additivity of~$\Phi$
follows from the inequalities\footnote{
Observe that
by assertion~(2) of Proposition~\ref{chain}
the gradient
$\nabla_h (u\circ\varphi)$
coincides with
$\nabla_h (\overline{u}_i\circ\varphi)$
a.e.\
on
$\varphi^{-1}(\operatorname{spt} \overline{u}_i)$,
for
$i=1, 2$,
and vanishes at almost all other points of~$\Omega$.}
\begin{multline*}
\sum\limits_{i=1}^2\Phi(V_i)-4\varepsilon\leq\sum\limits_{i=1}^2\int\limits_{\varphi^{-1}(V_i)}|\nabla_h(u_i\circ\varphi)(x)|^q\, dx-2\varepsilon \\
\leq \sum\limits_{i=1}^2\int\limits_{\varphi^{-1}(V_i)}|\nabla_h(\tilde{u}_i\circ\varphi)(x)|^q\, dx=\sum\limits_{i=1}^2\int\limits_{\varphi^{-1}(V_i)}|\nabla_h(\overline{u}_i\circ\varphi)(x)|^q\, dx\\
=\int\limits_{\varphi^{-1}(V_1\cup V_2)}|\nabla_h(u\circ\varphi)(x)|^q\, dx\leq \Phi(V_1\cup V_2).
\end{multline*}
\end{proof}

\begin{proof}[Proof of Theorem {\rm\ref{descr}}.]
The function
$\Phi\circ\varphi:U\mapsto \Phi(\varphi(U))$
of open sets
is quasi-additive
because $\Phi$
is quasi-additive by Lemma \ref{quasiaddit},
while $\varphi$
is a homeomorphism.
Lemma \ref{quasiaddit} also implies that
\eqref{Phii} holds for all
$u\in {\rm Lip}(Y)$.
Thus,
putting
$V=\varphi(B(z,r))$,
where
$B(z,r)\Subset\Omega$
is a~ball,
and dividing~\eqref{Phii} by
$|B(z,r)|$,
we find that 
$$
\frac{1}{|B(z,r)|}\int\limits_{B(z,r)}|\nabla_h (u\circ\varphi)(x)|^q\, dx\leq {\rm Lip}(u)^q\frac{\Phi(\varphi(B(z,r)))}{|B(z,r)|}
$$
for all
$u\in {\rm Lip}(Y)$.
Passing to the limit as
$r\to 0$,
by the Lebesgue differentiation theorem and Proposition~\ref{quasiadditprop}
we infer that
$$
|\nabla_h(u\circ\varphi)|\leq {\rm Lip}(u) ((\Phi\circ\varphi)')^{1/q}\quad
\text{a.e.,}
$$
where the majorant
$((\Phi\circ\varphi)')^{1/q}$
belongs to
$L_q(\Omega)$,
see Proposition~\ref{quasiadditprop}.
Hence,
$\varphi\in L^1_q(\Omega;Y)$
and
$|\nabla_0\varphi|^q\leq (\Phi\circ\varphi)^\prime$.
Proposition~\ref{quasiadditprop} yields
$$
\int\limits_{U}|\nabla_0\varphi|(x)^q\, dx\leq \int\limits_U (\Phi\circ\varphi)'(x)\, dx\leq \Phi(\varphi(U)).
$$

The reverse inequality is easier to establish.
If
${\rm Lip}(u)\leq 1$
and
${\rm dist}(\operatorname{spt} u,Y\setminus \varphi(U))>0$
then
$\nabla_h (u\circ\varphi)=0$
a.e.\ outside~$U$,
see assertion~(2) of Proposition~\ref{chain}.
Therefore,
$$
\int\limits_{\Omega}|\nabla_h(u\circ\varphi)(x)|^q\, dx=\int\limits_{U}|\nabla_h(u\circ\varphi)(x)|^q\, dx\leq \int\limits_{U}|\nabla_0\varphi|(x)^q\, dx.
$$
Since
$u:Y\to\mathbb{R}$
is arbitrary with
${\rm Lip}(u)\leq 1$
and
${\rm dist}(\operatorname{spt} u,Y\setminus \varphi(U))>0$,
this yields
$$
\Phi(\varphi(U))\leq \int\limits_{U}|\nabla_0\varphi|(x)^q\, dx.
$$
\end{proof}

Note the extreme simplicity of the case
$q=\infty$.

\begin{proposition}\label{q=oo}
Consider a~domain~$\Omega$
in some Carnot group~$\mathbb{G}$
and a~metric space
$(Y,d)$.
If a~mapping
$\varphi:\Omega\to Y$
induces the bounded composition operator
$$
\varphi^*:{\rm Lip}(Y)\to L^1_\infty(\Omega)
$$
then
$\varphi\in L^1_\infty(\Omega;Y)$
and
$\|\, |\nabla_0\varphi|\mid L_\infty(\Omega)\|=\|\varphi^*\|$.
\end{proposition}

\begin{proof}
By assumption,
for all
$u\in {\rm Lip}(Y)$
we have
$$
|\nabla_h (u\circ\varphi)(x)|\leq \|u\circ\varphi\mid L^1_\infty(\Omega)\|\leq {\rm Lip}(u)\cdot \|\varphi^*\|
$$
for a.e.
$x\in\Omega$.
This means that
$\varphi\in L^1_\infty(\Omega;Y)$
and
$\|\, |\nabla_0\varphi|\mid L_\infty(\Omega)\|\leq \|\varphi^*\|$.
The reverse inequality is established
by analogy with the argument made in the case
$q<\infty$
in the proof of Theorem~\ref{descr}.
\end{proof}

To close this section,
we apply standard arguments to verify that
the boundedness assumption for the composition operator
is superfluous.

\begin{definition}\label{Frechet}
A~real or complex linear space~$X$
is called an~$F$-{\it space}
whenever it is equipped with a~metric~$\rho$
such that:
\begin{itemize}
\item[(1)]
$\rho(x_1+x_3,x_2+x_3)=\rho(x_1,x_2)$
for all
$x_1$,
$x_2$,
$x_3\in X$;
\item[(2)]
$(X,\rho)$
is a~complete metric space;
\item[(3)]
the operations of addition and scalar multiplication
are continuous in the metric~$\rho$,
that is,
if a~sequence
$\{\alpha_n\}$
of numbers
converges to a~(finite) number~$\alpha$,
while two sequences
$\{x_n\}$
and
$\{y_n\}$
of vectors
converge in the~metric~$\rho$
to
$x\in X$
and
$y\in X$
respectively,
then
$$
\rho(\alpha_nx_n,\alpha x)\to 0,\quad \rho(x_n+y_n,x+y)\to 0.
$$
\end{itemize}
\end{definition}

\begin{proposition}\label{GraphClosed}
Consider a~domain~$\Omega$
in some Carnot group~$\mathbb{G}$
and a~metric space
$(Y,d)$.
Take
$1\leq q\leq \infty$.
If a~mapping
$\varphi:\Omega\to Y$
satisfies
$u\circ\varphi\in L^1_q(\Omega)$
for every
$u\in {\rm Lip}(Y)$
then the composition operator
$$
\varphi^*:{\rm Lip}(Y)\to L^1_q(\Omega)
$$
is bounded.
\end{proposition}

\begin{proof}
Take a point
$a\in Y$
and a monotone sequence
$\{\Omega_k\Subset \Omega\}$
of compact subdomains exhausting $\Omega$.
Consider the space
${\rm Lip}_a(Y)$
of all Lipschitz functions
$u:Y\to\mathbb{R}$
with
$u(a)=0$.
The functional
$u\mapsto{\rm Lip}(u)$
is a~norm on
${\rm Lip}_a(Y)$
making it a Banach space.
Furthermore,
define on
$L^1_q(\Omega)$
the structure of an $F$-space with the metric\footnote{
The convergence
$\rho(f_l,f)\to 0$
is equivalent to the collection of convergences
$f_l\to f$
in
$L_{1,{\rm loc}}(\Omega)$
and
$\nabla_h f_l\to \nabla_h f$
in
$L_q(\Omega)$.
This observation enables us to verify without effort
the validity of properties (2) and (3) in Definition \ref{Frechet}.}
$$
\rho(f,g)=\|f-g\mid L^1_q(\Omega)\|+\sum\limits_{k=1}^\infty\frac{1}{2^k}\frac{\|f-g\mid L_1(\Omega_k)\|}{1+\|f-g\mid L_1(\Omega_k)\|}.
$$

The composition operator
\begin{equation}\label{F-comp}
\varphi^*:\big({\rm Lip}_a(Y),{\rm Lip}(\cdot)\big)\to \big(L^1_q(\Omega),\rho\big)
\end{equation}
is closed.
Indeed,
take a~sequence
$\{u_l\in {\rm Lip}_a(Y)\}$
converging to some function~$u$
such that
the sequence
$\{\varphi^*u_l\}$
of images
converges to some function
$f\in L^1_q(\Omega)$
in the metric~$\rho$.
By local convergence in the mean,
upon passing to a~subsequence we may assume that
$\varphi^*u_l\to f$
a.e.
Since
$u_l(a)=0$
for all~$l$
and
${\rm Lip}(u_l-u)\to 0$,
we infer that
$u_l(y)\to u(y)$
for all
$y\in Y$,
and consequently
$u_l\circ\varphi\to u\circ\varphi$
everywhere in~$\Omega$.
This immediately implies that
$\varphi^*u=f$
a.e.

By the closed graph theorem for~$F$-spaces \cite[Theorem 2.15]{Rudin},
the operator \eqref{F-comp} is continuous.
Since the metric~$\rho$
dominates over
the
$L^1_q(\Omega)$-seminorm,
for this
$\varepsilon>0$
there exists
$\delta>0$
such that
$\| \varphi^*v\mid L^1_q(\Omega)\|\leq \varepsilon$
for
$v\in{\rm Lip}_a(Y)$
with
${\rm Lip}(v)\leq \delta$.
Inserting here
$v=\frac{\delta}{{\rm Lip}(u)}u$,
where
$0\neq u\in {\rm Lip}_a(Y)$,
and multiplying by
$\frac{{\rm Lip}(u)}{\delta}$,
we find that
for
$A=\frac{\varepsilon}{\delta}$
the inequality
\begin{equation}\label{a-bound}
\|u\circ\varphi\mid L^1_q(\Omega)\|\leq A\cdot{\rm Lip}(u)
\end{equation}
holds for all
$u\in {\rm Lip}_a(Y)$.
Since
$u=(u-u(a))+u(a)$,
while the seminorms of constant functions in
${\rm Lip}(Y)$
and
$L^1_q(\Omega)$
vanish,
\eqref{a-bound} holds for all
$u\in {\rm Lip}(Y)$.
\end{proof}

\section{Applications to the theory of composition operators of homogeneous Sobolev spaces}
\label{p<oo}
Given a~domain~$\Omega$
in some Carnot group~$\mathbb{G}$,
a~mapping
$\varphi\in\operatorname{ACL}(\Omega;\mathbb{G})$
is called a~mapping {\it with finite distortion}
whenever
$D_h\varphi=0$
a.e.\ on the zero set of the Jacobian
$Z=\{x\in\Omega\mid \det \widehat{D}\varphi(x)=0\}$.
For
$p\in[1;\infty)$
the {\it outer distortion function}
$K_{p}(\cdot,\varphi)$
is defined as
$$
K_{p}(x,\varphi) =
\begin{cases}
\dfrac{|D_h \varphi(x)|}{|\det \widehat{D}\varphi(x)|^{\frac{1}{p}}}
& \text{ if } \det \widehat{D} \varphi(x) \neq 0, \\
0
& \text{ otherwise}.
\end{cases}
$$

The space
$\operatorname{Lip}_{\operatorname{loc}}(\Omega)$
consists of all functions defined on~$\Omega$
and Lipschitz in the Carnot--Carath\'eodory metric
on every compact set
$K\Subset\Omega$.

The following theorem establishes a~close relation
between mappings with finite distortion and integrable distortion function
and the description of bounded composition operators
of homogeneous Sobolev spaces.

\begin{theorem}\label{descrCarnot}
A~homeomorphism
$\varphi:\Omega\to\Omega'$
of domains in an~arbitrary Carnot group~$\mathbb{G}$
induces the bounded composition operator
$$
\varphi^*:L^1_p(\Omega')\cap {\rm Lip}_{{\rm loc}}(\Omega')\to L^1_q(\Omega),
$$
where $1\leq q\leq p<\infty$,
as
$\varphi^*u=u\circ\varphi$
if and only if
\begin{itemize}
\item[\rm(1)]
$\varphi\in W^1_{q,{\rm loc}}(\Omega;\mathbb{G})$;
\item[\rm(2)]
$\varphi$
has finite distortion;
\item[\rm(3)]
$K_{p}(\cdot,\varphi)\in L_{\sigma}(\Omega)$,
where
$\frac{1}{\sigma}=\frac{1}{q}-\frac{1}{p}$,
as well as
$\sigma=\infty$
for
$q=p$.
\end{itemize}
Furthermore,
\begin{equation}\label{norms}
\|\varphi^*\|= \|K_p(\cdot,\varphi)\mid L_\sigma(\Omega)\|.
\end{equation}
\end{theorem}

\begin{remark}
Theorem~\ref{descrCarnot} is proved in
\cite[Theorem 2]{VodUkhlov1998}, \cite[Theorem 2]{VodEvs2022},
but instead of
\eqref{norms}
the cited articles only obtain the two-sided bounds
$$
\alpha_{q,p}\|K_p(\cdot,\varphi)\mid L_\sigma(\Omega)\|\leq \|\varphi^*\|\leq \|K_p(\cdot,\varphi)\mid L_\sigma(\Omega)\|,
$$
where
$\alpha_{q,p}\in (0;1)$
is a~constant.
Below we present a~new proof of the necessity of conditions (1)--(3),
which yields
\eqref{norms},
previously available \cite{Vod2022} only in the case
$\mathbb{G}=\mathbb{R}^n$.
It rests on Lemma~\ref{mainIneq}
and a~modification that lemma enables us to make
to the available approach of
\cite[Prop.~1]{VodUkhlov1998}, \cite{VodEvs2022}.
\end{remark}

\begin{proof}
Suppose that
$q<p$.
Under the conditions of Theorem~\ref{descrCarnot},
the function
\begin{multline*}
\Phi(V)=\bigg(\sup\bigg\{\frac{\|u\circ \varphi\mid L^1_q(\Omega)\|}{\|u\mid L^1_p(\Omega')\|}\, \Big\vert\, u\in L^1_p(\Omega')\cap {\rm Lip}_{\rm loc}(\Omega'),\\
u\vert_{\Omega'\setminus V}=0,\, \|u\mid L^1_p(\Omega')\|\neq 0 \bigg\}\bigg)^\sigma
\end{multline*}
of open sets
$V\subset\Omega'$
is quasi-additive
and
$\Phi(\Omega')= \|\varphi^*\|^\sigma$,
see \cite{VodUkhlov1998}, \cite[Lemma 3.1]{VodEvs2022}.
Here
$\frac{1}{\sigma}=\frac{1}{q}-\frac{1}{p}$.

Take some open set
$V\Subset\Omega'$.
Every~$1$-Lipschitz function
$u:\mathbb{G}\to\mathbb{R}$
belongs to
$L^1_\infty(\mathbb{G})$,
and
$|\nabla_h u|\leq 1$
a.e.
Therefore,
if
$u=0$
outside~$V$
then
$
\|u\mid L^1_p(\Omega')\|\leq |V|^{1/p}.
$
Hence,
by the definition of~$\Phi$
every~$1$-Lipschitz function
$u:\mathbb{G}\to \mathbb{R}$
with
${\rm dist}(\operatorname{spt} u,\mathbb{G}\setminus V)>0$
satisfies
\begin{equation}\label{B<2B-pq}
\left(\int\limits_{\varphi^{-1}(V)}|\nabla_h(u\circ\varphi)(x)|^q\, dx\right)^{1/q}\leq \Phi(V)^{1/\sigma}\|u\mid L^1_p(\Omega')\| \leq \Phi(V)^{1/\sigma}|V|^{1/p}.
\end{equation}
Lemma~\ref{mainIneq},
where we take
$(Y,d)=(\mathbb{G},d_{cc})$
and
$A=\Phi(V)^{q/\sigma}|V|^{q/p}$,
shows that
\eqref{B<2B-pq} holds for all~$1$-Lipschitz functions
$u:\mathbb{G}\to\mathbb{R}$.
Fix an~arbitrary~$1$-Lipschitz function
$u:\mathbb{G}\to\mathbb{R}$.
Insert into~\eqref{B<2B-pq}
the image
$V=\varphi(B(z,r))$
of a~sub-Riemannian ball
$B(z,r)\Subset \Omega$
and divide by
$|B(z,r)|^{1/q}$:
$$
\left(\frac{1}{|B(z,r)|}\int\limits_{B(z,r)}|\nabla_h(u\circ\varphi)(x)|^q\, dx\right)^{1/q}\leq \left(\frac{\Phi(\varphi(B(z,r)))}{|B(z,r)|}\right)^{1/\sigma}\left(\frac{|\varphi(B(z,r))|}{|B(z,r)|}\right)^{1/p}.
$$
Passing to the limit as
$r\to 0$,
by the Lebesgue differentiation theorem and Proposition~\ref{quasiadditprop}
we infer for every~$1$-Lipschitz function
$u:\mathbb{G}\to\mathbb{R}$
that
$$
|\nabla_h(u\circ\varphi)|\leq ((\Phi\circ\varphi)')^{1/\sigma}J(\cdot,\varphi)^{1/p}\quad \text{a.e.,}
$$
where
$\Phi\circ\varphi:U\mapsto \Phi(\varphi(U))$
is a~quasi-additive function,
$(\Phi\circ\varphi)'$
is its derivative,
while
$J(\cdot,\varphi)$
is the spatial derivative of the homeomorphism~$\varphi$,
see~\eqref{Jacob}.
The right-hand side of the last inequality lies in
$L_{q,{\rm loc}}(\Omega)$.
Indeed,
since
$q/p+q/\sigma=1$,
for
$U\Subset\Omega$
we deduce,
see assertion~(c) of Proposition~\ref{quasiadditprop},
that
\begin{multline*}
\int\limits_U((\Phi\circ\varphi)'(x))^{q/\sigma}J(x,\varphi)^{q/p}\, dx\leq \left(\int\limits_U(\Phi\circ\varphi)'(x)\, dx\right)^{q/\sigma}\left(\int\limits_U J(x,\varphi)\, dx\right)^{q/p}
\\
\leq (\Phi(\varphi(U))^{q/\sigma}|\varphi(U)|^{q/p}<\infty.
\end{multline*}
Hence,
$\varphi\in W^1_{q,{\rm loc}}(\Omega;\mathbb{G})$
and
\begin{equation}\label{final-major}
|\nabla_0\varphi|\leq ((\Phi\circ\varphi)')^{1/\sigma}J(\cdot,\varphi)^{1/p}\quad \text{a.e.}
\end{equation}
Since
$J(\cdot,\varphi)=|\det\widehat D\varphi|$
a.e.,
see subsection~\ref{prelim-quasiadd},
and\footnote{
The first inequality here
follows from the orthonormality of the vectors
$X_j$,
while the second one,
from the estimate
$|X_j\varphi|\leq |\nabla_0\varphi|$
in Proposition~\ref{R1p-ACL}.}
$$
|D_h\varphi|\leq \bigg(\sum_{j=1}^n |X_j\varphi|^2\bigg)^{1/2}\leq n^{1/2}|\nabla_0\varphi|\quad \text{a.e.},
$$
we infer that~$\varphi$
has finite distortion.
For mappings with finite distortion
$|\nabla_0\varphi|=|D_h\varphi|$
a.e.\ \cite[Lemma 3.4]{VodPav}.
Thus,
we can express~\eqref{final-major} as
$$
|D_h\varphi|\leq ((\Phi\circ\varphi)')^{1/\sigma}|\det\widehat D\varphi|^{1/p}.
$$
Hence,
$K_p(\cdot,\varphi)^\sigma\leq (\Phi\circ\varphi)'\in L_1(\Omega)$.
By Proposition~\ref{quasiadditprop},
every open set
$U\subset\Omega$
satisfies
$$
\int\limits_U K(x,\varphi)^\sigma(x)\, dx\leq \int\limits_U (\Phi\circ\varphi)'(x)\, dx\leq\Phi(\varphi(U))\leq \|\varphi^*\|^\sigma<\infty.
$$
The reverse inequalities
$\int\limits_{U}K_p(x,\varphi)^\sigma\, dx \geq \Phi(\varphi(U))$,
including the inequality
$$
\|K_{p}(\cdot,\varphi)\mid L_\sigma(\Omega)\|
\geq \Phi(\Omega')^{1/\sigma}= \|\varphi^*\|,
$$
are justified in \cite[Theorem 2]{VodUkhlov1998}.

In the case
$q=p$
the arguments simplify and avoid the quasi-additive function:
it suffices to replace
$\Phi(V)^{1/\sigma}$
by
$\|\varphi^*\|$
throughout the proof.
\end{proof}

\begin{corollary}\label{alpha=1}
Under the hypotheses of Theorem~{\rm\ref{descrCarnot}},
for
$q<p$
all open sets
$U\subset\Omega$
satisfy
$$
\int\limits_{U}K_p(x,\varphi)^\sigma\, dx = \Phi(\varphi(U)).
$$
In~particular,
$K_p(x,\varphi)^\sigma=(\Phi\circ\varphi)'(x)$
for a.e.\
$x\in\Omega$.
In addition,
$$
\|K_p(\cdot,\varphi)\mid L_\sigma(\Omega)\|=\|\varphi^*\|
$$
for all
$1\leq q\leq p<\infty$.
\end{corollary}

\section*{Acknowledgements}
S. V. Pavlov: The work is supported by the Mathematical Center in Akademgorodok
under the Agreement 075--15--2025--349
with the Ministry of Science and Higher Education of the Russian Federation.
\newline
S.~K.~Vodopyanov:
Working in the framework of the State Task
to the Sobolev Institute of Mathematics
from the Ministry of Higher Education and Science of the Russian Federation
(Project FWNF--2022--0006).

\end{document}